% E56, processed by citealice (2006-07-10) on Wed Jul 11 11:49:06 IDT 2007

\ifx\shlhetal\undefinedcontrolsequence\let\shlhetal\relax\fi
% To: shani@math.huji.ac.il
% Subject: e56
% CC: Saharon Shelah <shelah@math.huji.ac.il>, "Saharon Shelah's Office" <shlhetal@math.huji.ac.il>
% Date: Mon, 9 Jul 2007 09:25:38 -0400 (EDT)
% From: Alice Leonhardt <leonhard@math.rutgers.edu>
% Mime-Version: 1.0
% Content-Description: minor revision
% X-sliced-and-diced-by: 'savemail' 2.0, Oct 2005

\input amstex
% % \input mathdefs
%  *** start including mathdefs.tex *** 
\expandafter\ifx\csname mathdefs.tex\endcsname\relax
  \expandafter\gdef\csname mathdefs.tex\endcsname{}
\else \message{Hey!  Apparently you were trying to
  \string\input{mathdefs.tex} twice.   This does not make sense.} 
\errmessage{Please edit your file (probably \jobname.tex) and remove
any duplicate ``\string\input'' lines}\endinput\fi

%mathdefs.tex v1.3.2

%%% Changes from v1.0: footnote macros, warning for duplicated tags,
%%%   control sequences \( and \verbatimtags.
%%% From v1.2: \pretags, redefinition of \( using \ifinner, multi-part
%%%   equation numbering, control sequences \[, \references, and
%%%   \resetbracket. 
%%% From v1.3: \rm in \lastpart; write root of multi-part tag to .tgs 

%See file texdefs.doc for documentation.

\catcode`\X=12\catcode`\@=11

%Minor control sequences:
\def\n@wcount{\alloc@0\count\countdef\insc@unt}
\def\n@wwrite{\alloc@7\write\chardef\sixt@@n}
\def\n@wread{\alloc@6\read\chardef\sixt@@n}
\def\r@s@t{\relax}\def\v@idline{\par}\def\@mputate#1/{#1}
\def\l@c@l#1X{\firstpart.#1}\def\gl@b@l#1X{#1}\def\t@d@l#1X{{}}

%Creation of tag families and output of assignments and citations:
\def\crossrefs#1{\ifx\all#1\let\tr@ce=\all\else\def\tr@ce{#1,}\fi
   \n@wwrite\cit@tionsout\openout\cit@tionsout=\jobname.cit 
   \write\cit@tionsout{\tr@ce}\expandafter\setfl@gs\tr@ce,}
\def\setfl@gs#1,{\def\@{#1}\ifx\@\empty\let\next=\relax
   \else\let\next=\setfl@gs\expandafter\xdef
   \csname#1tr@cetrue\endcsname{}\fi\next}
\def\m@ketag#1#2{\expandafter\n@wcount\csname#2tagno\endcsname
     \csname#2tagno\endcsname=0\let\tail=\all\xdef\all{\tail#2,}
   \ifx#1\l@c@l\let\tail=\r@s@t\xdef\r@s@t{\csname#2tagno\endcsname=0\tail}\fi
   \expandafter\gdef\csname#2cite\endcsname##1{\expandafter
     \ifx\csname#2tag##1\endcsname\relax?\else\csname#2tag##1\endcsname\fi
     \expandafter\ifx\csname#2tr@cetrue\endcsname\relax\else
     \write\cit@tionsout{#2tag ##1 cited on page \folio.}\fi}
   \expandafter\gdef\csname#2page\endcsname##1{\expandafter
     \ifx\csname#2page##1\endcsname\relax?\else\csname#2page##1\endcsname\fi
     \expandafter\ifx\csname#2tr@cetrue\endcsname\relax\else
     \write\cit@tionsout{#2tag ##1 cited on page \folio.}\fi}
   \expandafter\gdef\csname#2tag\endcsname##1{\expandafter
      \ifx\csname#2check##1\endcsname\relax
      \expandafter\xdef\csname#2check##1\endcsname{}%
      \else\immediate\write16{Warning: #2tag ##1 used more than once.}\fi
      \multit@g{#1}{#2}##1/X%
      \write\t@gsout{#2tag ##1 assigned number \csname#2tag##1\endcsname\space
      on page \number\count0.}%
   \csname#2tag##1\endcsname}}

\def\multit@g#1#2#3/#4X{\def\t@mp{#4}\ifx\t@mp\empty%
      \global\advance\csname#2tagno\endcsname by 1 
      \expandafter\xdef\csname#2tag#3\endcsname
      {#1\number\csname#2tagno\endcsnameX}%
   \else\expandafter\ifx\csname#2last#3\endcsname\relax
      \expandafter\n@wcount\csname#2last#3\endcsname
      \global\advance\csname#2tagno\endcsname by 1 
      \expandafter\xdef\csname#2tag#3\endcsname
      {#1\number\csname#2tagno\endcsnameX}
      \write\t@gsout{#2tag #3 assigned number \csname#2tag#3\endcsname\space
      on page \number\count0.}\fi
   \global\advance\csname#2last#3\endcsname by 1
   \def\t@mp{\expandafter\xdef\csname#2tag#3/}%
   \expandafter\t@mp\@mputate#4\endcsname
   {\csname#2tag#3\endcsname\lastpart{\csname#2last#3\endcsname}}\fi}
\def\t@gs#1{\def\all{}\m@ketag#1e\m@ketag#1s\m@ketag\t@d@l p
\let\realscite\scite
\let\realstag\stag
   \m@ketag\gl@b@l r \n@wread\t@gsin
   \openin\t@gsin=\jobname.tgs \re@der \closein\t@gsin
   \n@wwrite\t@gsout\openout\t@gsout=\jobname.tgs }
\outer\def\localtags{\t@gs\l@c@l}
\outer\def\globaltags{\t@gs\gl@b@l}
\outer\def\newlocaltag#1{\m@ketag\l@c@l{#1}}
\outer\def\newglobaltag#1{\m@ketag\gl@b@l{#1}}

%Reading in tag information:
\newif\ifpr@ 
\def\m@kecs #1tag #2 assigned number #3 on page #4.%
   {\expandafter\gdef\csname#1tag#2\endcsname{#3}
   \expandafter\gdef\csname#1page#2\endcsname{#4}
   \ifpr@\expandafter\xdef\csname#1check#2\endcsname{}\fi}
\def\re@der{\ifeof\t@gsin\let\next=\relax\else
   \read\t@gsin to\t@gline\ifx\t@gline\v@idline\else
   \expandafter\m@kecs \t@gline\fi\let \next=\re@der\fi\next}
\def\pretags#1{\pr@true\pret@gs#1,,}
\def\pret@gs#1,{\def\@{#1}\ifx\@\empty\let\n@xtfile=\relax
   \else\let\n@xtfile=\pret@gs \openin\t@gsin=#1.tgs \message{#1} \re@der 
   \closein\t@gsin\fi \n@xtfile}

%Sections and subsections; local numbering:
\newcount\sectno\sectno=0\newcount\subsectno\subsectno=0
\newif\ifultr@local \def\ultralocal{\ultr@localtrue}
\def\firstpart{\number\sectno}
\def\lastpart#1{\ifcase#1 \or a\or b\or c\or d\or e\or f\or g\or h\or 
   i\or k\or l\or m\or n\or o\or p\or q\or r\or s\or t\or u\or v\or w\or 
   x\or y\or z \fi}

\def\resetall{\global\advance\sectno by 1\subsectno=0
   \gdef\firstpart{\number\sectno}\r@s@t}
\def\resetsub{\global\advance\subsectno by 1
   \gdef\firstpart{\number\sectno.\number\subsectno}\r@s@t}
\def\newsection#1\par{\resetall\vskip0pt plus.3\vsize\penalty-250
   \vskip0pt plus-.3\vsize\bigskip\bigskip
   \message{#1}\leftline{\bf#1}\nobreak\bigskip}
\def\subsection#1\par{\ifultr@local\resetsub\fi
   \vskip0pt plus.2\vsize\penalty-250\vskip0pt plus-.2\vsize
   \bigskip\smallskip\message{#1}\leftline{\bf#1}\nobreak\medskip}

%jj tags:
% On Andrzej's request:  we want to be able 
% to show tags as in noverbatim, with verbatim in the margin,
% and cites as in verbatim, with nonverbatim in the margin
% mg -- July 2000

\newdimen\marginshift

\newdimen\margindelta
\newdimen\marginmax
\newdimen\marginmin

\def\margininit{       
\marginmax=3 true cm                  % how much room, approximately
				      
\margindelta=0.1 true cm              % distance between entries
\marginmin=0.1true cm                 % where will leftmost entry be
\marginshift=\marginmin
}    % we cannot execute this right now, since 
     % there may be a \magnification coming later in the 
     % main file.   So we call \margininit at the end of 
     % alice2jlem

\def\t@gsjj#1,{\def\@{#1}\ifx\@\empty\let\next=\relax\else\let\next=\t@gsjj
   \def\@@{p}\ifx\@\@@\else
   \expandafter\gdef\csname#1cite\endcsname##1{\citejj{##1}}
   \expandafter\gdef\csname#1page\endcsname##1{?}
   \expandafter\gdef\csname#1tag\endcsname##1{\tagjj{##1}}\fi\fi\next}
\newif\ifshowstuffinmargin
\showstuffinmarginfalse
\def\jjtags{\ifx\shlhetal\relax 
      % so this is a public version --> no-op 
  \else
\ifx\shlhetal\undefinedcontrolseq
      % again, this is a public version --> no-op 
\else
\showstuffinmargintrue
\ifx\all\relax\else\expandafter\t@gsjj\all,\fi\fi \fi
}

\def\tagjj#1{\realstag{#1}\oldmginpar{\zeigen{#1}}}
\def\citejj#1{\rechnen{#1}\mginpar{\zeigen{#1}}}     % modified Sep 02, saharon's suggestion

\def\rechnen#1{\expandafter\ifx\csname stag#1\endcsname\relax ??\else
                           \csname stag#1\endcsname\fi}

\newdimen\theight

\def\marginfont{\sevenrm}

\def\trymarginbox#1{\setbox0=\hbox{\marginfont\hskip\marginshift #1}%
		\global\marginshift\wd0 
		\global\advance\marginshift\margindelta}

\def \oldmginpar#1{%
\ifvmode\setbox0\hbox to \hsize{\hfill\rlap{\marginfont\quad#1}}%
\ht0 0cm
\dp0 0cm
\box0\vskip-\baselineskip
\else 
             \vadjust{\trymarginbox{#1}%
		\ifdim\marginshift>\marginmax \global\marginshift\marginmin
			\trymarginbox{#1}%
                \fi
             \theight=\ht0
             \advance\theight by \dp0    \advance\theight by \lineskip
             \kern -\theight \vbox to \theight{\rightline{\rlap{\box0}}%
\vss}}\fi}

\newdimen\upordown
\global\upordown=8pt
\font\tinyfont=cmtt8 % scaled 700
\def\mginpar#1{\smash{\hbox to 0cm{\kern-10pt\raise7pt\hbox{\tinyfont #1}\hss}}}
% testing, october 2005, mg
\def\mginpar#1{{\hbox to 0cm{\kern-10pt\raise\upordown\hbox{\tinyfont #1}\hss}}\global\upordown-\upordown}

% \def\mginpar#1{mg-#1-mg }

%Verbatim tags:
\def\t@gsoff#1,{\def\@{#1}\ifx\@\empty\let\next=\relax\else\let\next=\t@gsoff
   \def\@@{p}\ifx\@\@@\else
   \expandafter\gdef\csname#1cite\endcsname##1{\zeigen{##1}}
   \expandafter\gdef\csname#1page\endcsname##1{?}
   \expandafter\gdef\csname#1tag\endcsname##1{\zeigen{##1}}\fi\fi\next}
\def\verbatimtags{\showstuffinmarginfalse
\ifx\all\relax\else\expandafter\t@gsoff\all,\fi}
\def\zeigen#1{\hbox{$\scriptstyle\langle$}#1\hbox{$\scriptstyle\rangle$}}

% % \def\margincite#1{\ifshowstuffinmargin\mginpar{\rechnen{#1}}\fi}
%  changed, april 2003, mg: we now have always the verbatim tag in the margin

%Equation numbering:
\def\(#1){\edef\dot@g{\ifmmode\ifinner(\hbox{\noexpand\etag{#1}})
   \else\noexpand\eqno(\hbox{\noexpand\etag{#1}})\fi
   \else(\noexpand\ecite{#1})\fi}\dot@g}

%Reference numbering:
\newif\ifbr@ck
\def\eat#1{}
\def\[#1]{\br@cktrue[\br@cket#1'X]}
\def\br@cket#1'#2X{\def\temp{#2}\ifx\temp\empty\let\next\eat
   \else\let\next\br@cket\fi
   \ifbr@ck\br@ckfalse\br@ck@t#1,X\else\br@cktrue#1\fi\next#2X}
\def\br@ck@t#1,#2X{\def\temp{#2}\ifx\temp\empty\let\neext\eat
   \else\let\neext\br@ck@t\def\temp{,}\fi
   \def\teemp{#1}\ifx\teemp\empty\else\rcite{#1}\fi\temp\neext#2X}
\def\resetbr@cket{\gdef\[##1]{[\rtag{##1}]}}
\def\references{\resetbr@cket\newsection References\par}

%Footnotes:
\newtoks\symb@ls\newtoks\s@mb@ls\newtoks\p@gelist\n@wcount\ftn@mber
    \ftn@mber=1\newif\ifftn@mbers\ftn@mbersfalse\newif\ifbyp@ge\byp@gefalse
\def\defm@rk{\ifftn@mbers\n@mberm@rk\else\symb@lm@rk\fi}
\def\n@mberm@rk{\xdef\m@rk{{\the\ftn@mber}}%
    \global\advance\ftn@mber by 1 }
\def\rot@te#1{\let\temp=#1\global#1=\expandafter\r@t@te\the\temp,X}
\def\r@t@te#1,#2X{{#2#1}\xdef\m@rk{{#1}}}
\def\b@@st#1{{$^{#1}$}}\def\str@p#1{#1}
\def\symb@lm@rk{\ifbyp@ge\rot@te\p@gelist\ifnum\expandafter\str@p\m@rk=1 
    \s@mb@ls=\symb@ls\fi\write\f@nsout{\number\count0}\fi \rot@te\s@mb@ls}
\def\byp@ge{\byp@getrue\n@wwrite\f@nsin\openin\f@nsin=\jobname.fns 
    \n@wcount\currentp@ge\currentp@ge=0\p@gelist={0}
    \re@dfns\closein\f@nsin\rot@te\p@gelist
    \n@wread\f@nsout\openout\f@nsout=\jobname.fns }
\def\m@kelist#1X#2{{#1,#2}}
\def\re@dfns{\ifeof\f@nsin\let\next=\relax\else\read\f@nsin to \f@nline
    \ifx\f@nline\v@idline\else\let\t@mplist=\p@gelist
    \ifnum\currentp@ge=\f@nline
    \global\p@gelist=\expandafter\m@kelist\the\t@mplistX0
    \else\currentp@ge=\f@nline
    \global\p@gelist=\expandafter\m@kelist\the\t@mplistX1\fi\fi
    \let\next=\re@dfns\fi\next}
\def\symbols#1{\symb@ls={#1}\s@mb@ls=\symb@ls} 
\def\bigsymbol{\textstyle}
\symbols{\bigsymbol\ast,\dagger,\ddagger,\sharp,\flat,\natural,\star}
\def\ftnumbers{\ftn@mberstrue} \def\ftsymbols{\ftn@mbersfalse}
\def\paginal{\byp@ge} \def\resetftnumbers{\ftn@mber=1}
\def\ftnote#1{\defm@rk\expandafter\expandafter\expandafter\footnote
    \expandafter\b@@st\m@rk{#1}}

%Miscellaneous macros:
\long\def\jump#1\endjump{}
\def\ssum{\mathop{\lower .1em\hbox{$\textstyle\Sigma$}}\nolimits}

\def\qed{\nobreak\kern 1em \vrule height .5em width .5em depth 0em}
\def\newneq{\hbox{\rlap{\hbox to 1\wd9{\hss$=$\hss}}\raise .1em 
   \hbox to 1\wd9{\hss$\scriptscriptstyle/$\hss}}}
\def\subsetne{\setbox9 = \hbox{$\subset$}\mathrel{\hbox{\rlap
   {\lower .4em \newneq}\raise .13em \hbox{$\subset$}}}}
\def\supsetne{\setbox9 = \hbox{$\subset$}\mathrel{\hbox{\rlap
   {\lower .4em \newneq}\raise .13em \hbox{$\supset$}}}}

%Blackboard bold:
\def\vbar{\mathchoice{\vrule height6.3ptdepth-.5ptwidth.8pt\kern-.8pt}
   {\vrule height6.3ptdepth-.5ptwidth.8pt\kern-.8pt}
   {\vrule height4.1ptdepth-.35ptwidth.6pt\kern-.6pt}
   {\vrule height3.1ptdepth-.25ptwidth.5pt\kern-.5pt}}
\def\f@dge{\mathchoice{}{}{\mkern.5mu}{\mkern.8mu}}
\def\b@c#1#2{{\rm \mkern#2mu\vbar\mkern-#2mu#1}}
\def\b@b#1{{\rm I\mkern-3.5mu #1}}
\def\b@a#1#2{{\rm #1\mkern-#2mu\f@dge #1}}
\def\bb#1{{\count4=`#1 \advance\count4by-64 \ifcase\count4\or\b@a A{11.5}\or
   \b@b B\or\b@c C{5}\or\b@b D\or\b@b E\or\b@b F \or\b@c G{5}\or\b@b H\or
   \b@b I\or\b@c J{3}\or\b@b K\or\b@b L \or\b@b M\or\b@b N\or\b@c O{5} \or
   \b@b P\or\b@c Q{5}\or\b@b R\or\b@a S{8}\or\b@a T{10.5}\or\b@c U{5}\or
   \b@a V{12}\or\b@a W{16.5}\or\b@a X{11}\or\b@a Y{11.7}\or\b@a Z{7.5}\fi}}

\catcode`\X=11 \catcode`\@=12

% Sep 2003, mg:

% % \input 300stuff
%  *** start including 300stuff.tex *** 

% definitions needed to process the 300x papers  (plus 88 etc) 

\let\thischap\jobname

\def\partof#1{\csname returnthe#1part\endcsname}
\def\chapof#1{\csname returnthe#1chap\endcsname}

\def\setchapter#1,#2,#3;{% 
  \expandafter\def\csname returnthe#1part\endcsname{#2}%
  \expandafter\def\csname returnthe#1chap\endcsname{#3}%
}

\setchapter 300a,A,II.A;
\setchapter 300b,A,II.B;
\setchapter 300c,A,II.C;
\setchapter 300d,A,II.D;
\setchapter 300e,A,II.E;
\setchapter 300f,A,II.F;
\setchapter 300g,A,II.G;
\setchapter  E53,B,N;
\setchapter  88r,B,I;
\setchapter  600,B,III;
\setchapter  705,B,IV;
\setchapter  734,B,V;

\def\cprefix#1{%     \cprefix{300b} generates II.  or "nothing"
 %\leavevmode\vrule width 2cm height 1cm depth 0cm
\edef\theotherpart{\partof{#1}}\edef\theotherchap{\chapof{#1}}%
\ifx\theotherpart\thispart
   \ifx\theotherchap\thischap % nothing
    \else % same part, different chap
     \theotherchap%
    \fi
   \else % different part
     \theotherchap\fi}

 % \sectioncite[\S4]{300b}  ->  II.\S4   or A.II.\S4
\def\sectioncite[#1]#2{%
     \cprefix{#2}#1}

% now define \thispart
\edef\thispart{\partof{\thischap}}
\edef\thischap{\chapof{\thischap}}

\def\lastpage of '#1' is #2.{\expandafter\def\csname lastpage#1\endcsname{#2}}

%  *** end including 300stuff.tex *** 

\def\spuriousreset{}

%  *** end including mathdefs.tex *** 
% % \input citeadd
%  *** start including citeadd.tex *** 
%   citeadd -- a few additions for 
% files from alice that were procesed with "citealice"

\expandafter\ifx\csname citeadd.tex\endcsname\relax
\expandafter\gdef\csname citeadd.tex\endcsname{}
\else \message{Hey!  Apparently you were trying to
\string\input{citeadd.tex} twice.   This does not make sense.} 
\errmessage{Please edit your file (probably \jobname.tex) and remove
any duplicate ``\string\input'' lines}\endinput\fi

%  *** end including citeadd.tex *** 
\sectno=-1   % start with sect 0
\localtags
\jjtags
\NoBlackBoxes
\define\mr{\medskip\roster}
\define\sn{\smallskip\noindent}

\define\ub{\underbar}

\define\ermn{\endroster\medskip\noindent}

\define \nl{\newline}
\magnification=\magstep 1
\documentstyle{amsppt}
% % \input alice2000
%  *** start including alice2000.tex *** 
% This file should be inputted whenever we use amsppt.sty and 
% old tex.  
%  Here we redefine \subjclass (use 1991 instead of 2000, otherwise 
% the following definition comes directly from 
%% 
%%              `amsppt.sty', generated 
%% on <1997/2/2> with the docstrip utility (2.2i).
%% 
%% The original source files were:
%% 
%% amsppt.doc 
%%% ====================================================================
%%% @AMSTeX-style-file{
%%%   filename  = "amsppt.sty",
%%%   version   = "2.1h",
%%%   date      = "1997/02/02",
%%%   time      = "09:27:44 EST",
%%%   checksum  = "56844 3264 16617 137829",
%%%   author    = "American Mathematical Society",
%%%   address   = "PO Box 6248, Providence, RI 02940-6248, USA",
%%%   telephone = "401-455-4080 or (in the USA) 800-321-4AMS",

{    % the braces make the catcode-change local. 
\catcode`@11

\ifx\alicetwothousandloaded@\relax
  \endinput\else\global\let\alicetwothousandloaded@\relax\fi

\gdef\subjclass{\let\savedef@\subjclass
 \def\subjclass##1\endsubjclass{\let\subjclass\savedef@
   \toks@{\def\usualspace{{\rm\enspace}}\eightpoint}%
   \toks@@{##1\unskip.}%
   \edef\thesubjclass@{\the\toks@
     \frills@{{\noexpand\rm2000 {\noexpand\it Mathematics Subject
       Classification}.\noexpand\enspace}}%
     \the\toks@@}}%
  \nofrillscheck\subjclass}
} 

%  *** end including alice2000.tex *** 
% % \input alice2jlem
%  *** start including alice2jlem.tex *** 
%% # Keywords  Input file to be used for texing Alice's files

\expandafter\ifx\csname alice2jlem.tex\endcsname\relax
  \expandafter\xdef\csname alice2jlem.tex\endcsname{\the\catcode`@}
\else \message{Hey!  Apparently you were trying to
\string\input{alice2jlem.tex}  twice.   This does not make sense.}
\errmessage{Please edit your file (probably \jobname.tex) and remove
any duplicate ``\string\input'' lines}\endinput\fi

% % \input bib4plain
%  *** start including bib4plain.tex *** 
\expandafter\ifx\csname bib4plain.tex\endcsname\relax
  \expandafter\gdef\csname bib4plain.tex\endcsname{}
\else \message{Hey!  Apparently you were trying to \string\input
  bib4plain.tex twice.   This does not make sense.}
\errmessage{Please edit your file (probably \jobname.tex) and remove
any duplicate ``\string\input'' lines}\endinput\fi

%  This file should be inputted if you want to use 
%  bibtex fom within plain TeX. 
      % Not really need for standard
       % bibtex files, but these commands
\def\renewcommand{\newcommand}	       % are used in our literal-unsrt.bst
\edef\cite{\the\catcode`@}%
\catcode`@ = 11
\let\@oldatcatcode = \cite
\chardef\@letter = 11
\chardef\@other = 12
%
%
% Next come some things that will be useful later.
%
% Make an outer definition into an inner one (due to Chris Thompson).
% The arguments should be the control sequence to be defined, and the
% new of the \outer control sequence, as characters; the control
% sequence #1 is defined to be just the same as \csname#2\endcsname, but
% not \outer.  For example, \@innerdef\innernewcount{newcount} would
% define \innernewcount to be a non-outer version of \newcount.
%
\def\@innerdef#1#2{\edef#1{\expandafter\noexpand\csname #2\endcsname}}%
%
% We use \@innerdef to make some of our allocations local, because
% Eplain includes our code inside a conditional.  We put @'s in the
% names to minimize the (already small) chance of conflicts.
%
\@innerdef\@innernewcount{newcount}%
\@innerdef\@innernewdimen{newdimen}%
\@innerdef\@innernewif{newif}%
\@innerdef\@innernewwrite{newwrite}%
%
%
% Swallow one parameter.
%
\def\@gobble#1{}%
%
%
% Use TeX 3.0's \inputlineno to get the line number, for better error
% messages, but if we're using an old version of TeX, don't do anything.
%
\ifx\inputlineno\@undefined
   \let\@linenumber = \empty % Pre-3.0.
\else
   \def\@linenumber{\the\inputlineno:\space}%
\fi
%
%
% The following macro \@futurenonspacelet (from the TeXbook) behaves
% essentially like \futurelet except that it discards any implicit or
% explicit space tokens that intervene before a nonspace is scanned:
%
\def\@futurenonspacelet#1{\def\cs{#1}%
   \afterassignment\@stepone\let\@nexttoken=
}%
\begingroup % The grouping here avoids stepping on an outside use of `\\'.
\def\\{\global\let\@stoken= }%
\\ % now \@stoken is a space token (\\ is a control symbol, so that
   % space after it is seen).
\endgroup
\def\@stepone{\expandafter\futurelet\cs\@steptwo}%
\def\@steptwo{\expandafter\ifx\cs\@stoken\let\@@next=\@stepthree
   \else\let\@@next=\@nexttoken\fi \@@next}%
\def\@stepthree{\afterassignment\@stepone\let\@@next= }%
%
%
% \@getoptionalarg\CS gets an optional argument from the input, enclosed
% in brackets, then expands \CS.  We set \@optionalarg to \empty if we
% don't find one, otherwise to the text of the argument.  This assumes
% the brackets don't have a funny category code.
%
\def\@getoptionalarg#1{%
   \let\@optionaltemp = #1%
   \let\@optionalnext = \relax
   \@futurenonspacelet\@optionalnext\@bracketcheck
}%
%
% The \expandafter's in this macro let us avoid the use of \aftergroup,
% which is somewhat more expensive.
%
\def\@bracketcheck{%
   \ifx [\@optionalnext
      \expandafter\@@getoptionalarg
   \else
      \let\@optionalarg = \empty
      % We can't do the \temp after the \fi, because then the \temp gets
      % in the way of reading the optional argument from the input, if
      % we do have one.
      \expandafter\@optionaltemp
   \fi
}%
\def\@@getoptionalarg[#1]{%
   \def\@optionalarg{#1}%
   \@optionaltemp
}%
%
%
% From LaTeX.
%
\def\@nnil{\@nil}%
\def\@fornoop#1\@@#2#3{}%
\def\@for#1:=#2\do#3{%
   \edef\@fortmp{#2}%
   \ifx\@fortmp\empty \else
      \expandafter\@forloop#2,\@nil,\@nil\@@#1{#3}%
   \fi
}%
\def\@forloop#1,#2,#3\@@#4#5{\def#4{#1}\ifx #4\@nnil \else
       #5\def#4{#2}\ifx #4\@nnil \else#5\@iforloop #3\@@#4{#5}\fi\fi
}%
\def\@iforloop#1,#2\@@#3#4{\def#3{#1}\ifx #3\@nnil
       \let\@nextwhile=\@fornoop \else
      #4\relax\let\@nextwhile=\@iforloop\fi\@nextwhile#2\@@#3{#4}%
}%
%
%
% This macro tests if a file \jobname.#1 exists, and sets \if@fileexists
% appropriately.  If an optional argument is given, it is used as the
% root part of the filename instead of \jobname.
%
\@innernewif\if@fileexists
\def\@testfileexistence{\@getoptionalarg\@finishtestfileexistence}%
\def\@finishtestfileexistence#1{%
   \begingroup
      \def\extension{#1}%
      \immediate\openin0 =
         \ifx\@optionalarg\empty\jobname\else\@optionalarg\fi
         \ifx\extension\empty \else .#1\fi
         \space
      \ifeof 0
         \global\@fileexistsfalse
      \else
         \global\@fileexiststrue
      \fi
      \immediate\closein0
   \endgroup
}%
%
%
%% [[[start of BibTeX-specific stuff]]]
%
% Now come the four main LaTeX commands and their associated .aux
% commands.  Just as in LaTeX, \bibliographystyle defines the BibTeX
% style name (.bst file, that is), and \bibliography defines the
% database (.bib) file(s).  The corresponding .aux-file commands are
% \bibstyle and \bibdata, which are there only for BibTeX's (but not
% LaTeX's) use.
%
\def\bibliographystyle#1{%
   \@readauxfile
   \@writeaux{\string\bibstyle{#1}}%
}%
\let\bibstyle = \@gobble
%
% As well as writing the \bibdata command to tell BibTeX which .bib
% files to read, we read the .bbl file that BibTeX (or a person,
% conceivably) has produced.  We use \bblfilebasename as the root of the
% filename to read; this defaults to \jobname.
%
\let\bblfilebasename = \jobname
\def\bibliography#1{%
   \@readauxfile
   \@writeaux{\string\bibdata{#1}}%
   \@testfileexistence[\bblfilebasename]{bbl}%
   \if@fileexists
      % We just output a non-discardable item (the `whatsit' with the
      % \bibdata command).  This means that the glue that will be
      % inserted next (\parskip or \baselineskip, most likely) will be a
      % legal breakpoint.  Most likely, this is after some kind of
      % heading, where we don't want to allow a page break.  So:
      \nobreak
      \@readbblfile
   \fi
}%
\let\bibdata = \@gobble
%
% The \nocite{label,label,...} command writes its argument to \@auxfile,
% unless instructed not to, but produces no text in the document.  Both
% the \nocite and \cite commands produce \citation commands in the .aux file.
%
\def\nocite#1{%
   \@readauxfile
   \@writeaux{\string\citation{#1}}%
}%
\@innernewif\if@notfirstcitation
%
% \cite[note]{label,label,...} produces the citations for the labels as
% well.  If the optional argument `note' is present, it's added after
% the labels.  Since \cite calls \nocite to do its .aux-file writing,
% \cite doesn't need to call \@readauxfile (\nocite does).
%
\def\cite{\@getoptionalarg\@cite}%
%
% Typeset the citations for the labels in #1, followed by the note, if
% it exists.  To change the citation's format in the text, redefine one
% or more `\print...' macros, whose defaults appear later in this file.
%
\def\@cite#1{%
   % Remember the optional argument, in case one of the macros we call
   % below ends up looking for an optional argument itself.  For
   % example, if a \cite[note] triggers reading the .aux file, then the
   % [note] would be clobbered, since \@testfileexistence looks for an
   % optional arg.
   \let\@citenotetext = \@optionalarg
   % Start printing the text, beginning with a left bracket by default.
   \printcitestart
   % It's complicated, but because \nocite puts a `whatsit' onto the list,
   % \nocite should follow \printcitestart.  It's conceivable, but very
   % unlikely, that this `whatsit' will cause a problem (glue that doesn't
   % disappear when you want it to is the most likely symptom), requiring
   % a change either to \printcitestart or to the label that the .bst file
   % produces.
   \nocite{#1}%
   \@notfirstcitationfalse
   \@for \@citation :=#1\do
   {%
      \expandafter\@onecitation\@citation\@@
   }%
   \ifx\empty\@citenotetext\else
      \printcitenote{\@citenotetext}%
   \fi
   \printcitefinish
}%
\newif\ifweareinprivate
\weareinprivatetrue
\ifx\shlhetal\undefinedcontrolseq\weareinprivatefalse\fi
\ifx\shlhetal\relax\weareinprivatefalse\fi
\def\@onecitation#1\@@{%
   \if@notfirstcitation
      \printbetweencitations
   \fi
   \expandafter \ifx \csname\@citelabel{#1}\endcsname \relax
      \if@citewarning
         \message{\@linenumber Undefined citation `#1'.}%
      \fi
      % Give it a dummy definition:
     \ifweareinprivate
      \expandafter\gdef\csname\@citelabel{#1}\endcsname{%
% Change: marginal remark added, goldstrn@math.huji.ac.il, 
% goldstern@tuwien.ac.at, May 1996 mg
%  !!! change !!!
\strut 
\vadjust{\vskip-\dp\strutbox
\vbox to 0pt{\vss\parindent0cm \leftskip=\hsize 
\advance\leftskip3mm
\advance\hsize 4cm\strut\openup-4pt 
\rightskip 0cm plus 1cm minus 0.5cm ?  #1 ?\strut}}
         {\tt
            \escapechar = -1
            \nobreak\hskip0pt\pfeilsw%\special{ps:100 0 0 setrgbcolor }%
            \expandafter\string\csname#1\endcsname
                     %\special{ps:0 0 0 setrgbcolor }
             \pfeilso
            \nobreak\hskip0pt
         }%
      }%
     \else  % ifweareinprivate = false
      \expandafter\gdef\csname\@citelabel{#1}\endcsname{%
            {\tt\expandafter\string\csname#1\endcsname}
      }%
     \fi  % ifweareinprivate
   \fi
   % Now produce the text, whether it was undefined or not.
   \csname\@citelabel{#1}\endcsname
   \@notfirstcitationtrue
}%
%
% Given a label `foo', the macro `\b@foo' is supposed to
% hold the text that should be produced.
%
\def\@citelabel#1{b@#1}%
%
% So, how does a citation label get defined?  When we read the .bbl file
% (below), a \bibitem writes out a \@citedef command.  And when we read
% the \@citedef, we define \@citelabel{#1}, where #1 is the user's
% label.
%
\def\@citedef#1#2{\expandafter\gdef\csname\@citelabel{#1}\endcsname{#2}}%
%
%
% Reading the .bbl file also produces the typeset bibliography.  Please
% notice, however, that we do not produce the title for the references
% (e.g., `References'), as LaTeX does.  The formatting and spacing of
% that title, whether it should go into the headline, and so on, are all
% things determined by your format.  We cannot know those things in
% advance.  If you wish, you can define \bblhook to produce the title.
% Or just do it before the \bibliography command.
%
\def\@readbblfile{%
   % Define a counter to tell us which item number we are on, unless
   % we've already defined it (because the document has more than one
   % bibliography).
   \ifx\@itemnum\@undefined
      \@innernewcount\@itemnum
   \fi
   \begingroup
      \def\begin##1##2{%
         % ##1 is just `thebibliography'.
         % ##2 is the widest label.
         % We set (new dimen) \biblabelwidth based on the widest label
         \setbox0 = \hbox{\biblabelcontents{##2}}%
         \biblabelwidth = \wd0
      }%
      \def\end##1{}% ##1 is `thebibliography' again.
      %
      % Here we have two possibilities:
      % \bibitem[typesetlabel]{citationlabel}
      % \bibitem{citationlabel}
      % If we have the second of these, the citations are numbered, starting
      % from one; we use our own count register \@itemnum for this.
      %
      \@itemnum = 0
      \def\bibitem{\@getoptionalarg\@bibitem}%
      \def\@bibitem{%
         \ifx\@optionalarg\empty
            \expandafter\@numberedbibitem
         \else
            \expandafter\@alphabibitem
         \fi
      }%
      \def\@alphabibitem##1{%
         % Need \xdef here for various reasons.
         \expandafter \xdef\csname\@citelabel{##1}\endcsname {\@optionalarg}%
         % Left-justify alpha labels, unless \biblabel{pre,post}contents
         % are already defined.
         \ifx\biblabelprecontents\@undefined
            \let\biblabelprecontents = \relax
         \fi
         \ifx\biblabelpostcontents\@undefined
            \let\biblabelpostcontents = \hss
         \fi
         \@finishbibitem{##1}%
      }%
      \def\@numberedbibitem##1{%
         \advance\@itemnum by 1
         \expandafter \xdef\csname\@citelabel{##1}\endcsname{\number\@itemnum}%
         % Right-justify numeric labels, unless \biblabel{pre,post}contents
         % are already defined.
         \ifx\biblabelprecontents\@undefined
            \let\biblabelprecontents = \hss
         \fi
         \ifx\biblabelpostcontents\@undefined
            \let\biblabelpostcontents = \relax
         \fi
         \@finishbibitem{##1}%
      }%
      \def\@finishbibitem##1{%
         \biblabelprint{\csname\@citelabel{##1}\endcsname}%
         \@writeaux{\string\@citedef{##1}{\csname\@citelabel{##1}\endcsname}}%
         \ignorespaces
      }%
      %
      % Do the printing (we're producing the bibliography, remember).
      %
      \let\em = \bblem
      \let\newblock = \bblnewblock
      \let\sc = \bblsc
      % Punctuation won't affect spacing;
      \frenchspacing
      % the penalties below are from LaTeX's [article,book,report].sty;
      \clubpenalty = 4000 \widowpenalty = 4000
      % the next two values come from LaTeX's \sloppy command;
      \tolerance = 10000 \hfuzz = .5pt
      \everypar = {\hangindent = \biblabelwidth
                      \advance\hangindent by \biblabelextraspace}%
      \bblrm
      % the \parskip is a guess at what looks good;
      \parskip = 1.5ex plus .5ex minus .5ex
      % and the space between label and text comes from LaTeX's \labelsep.
      \biblabelextraspace = .5em
      \bblhook
      \input \bblfilebasename.bbl
   \endgroup
}%
%
% The widest label's width is useful for redefining \biblabelprint;
% you redefine \biblabelwidth, in effect, by redefining the
% \biblabelcontents macro that appears below.  And \biblabelextraspace,
% which is redefinable inside \bblhook, is added to \biblabelwidth to
% determine the amount of hanging indentation.
%
\@innernewdimen\biblabelwidth
\@innernewdimen\biblabelextraspace
%
% Now come the main macros that are related to the printing of the
% bibliography.  Since you might want to redefine them, they are given
% default definitions outside of \@readbblfile.
%
% The first one controls the printing of a bibliography entry's label.
% If you change it, make sure that it starts with something like
% \noindent or \indent or \leavevmode that puts TeX into horizontal mode
% (even if the label itself is empty); otherwise, the hanging
% indentation will get messed up in certain circumstances.
%
\def\biblabelprint#1{%
   \noindent
   \hbox to \biblabelwidth{%
      \biblabelprecontents
      \biblabelcontents{#1}%
      \biblabelpostcontents
   }%
   \kern\biblabelextraspace
}%
%
% If you are using numeric labels, and you want them left-justified
% (numeric labels by default are right-justified), do something like:
%     \def\biblabelprecontents{\relax}
%     \def\biblabelpostcontents{\hss}
%
% By default the labels are typeset in \bblrm, and enclosed in brackets.
%
\def\biblabelcontents#1{{\bblrm [#1]}}%
%
% The main text, too, is typeset using \bblrm, which is \rm by default.
%
\def\bblrm{\rm}%
%
% Emphasis for producing, e.g., titles, is done with \it by default.
%
\def\bblem{\it}%
%
% Some styles use a caps-and-small-caps font for author names.  LaTeX
% defines an \sc command but plain TeX doesn't, so we need one here.
% The definition below doesn't load the font unless it's needed, but it
% tries to load only the 10pt version, because it might not exist at
% other point sizes.
%
\def\bblsc{\ifx\@scfont\@undefined
              \font\@scfont = cmcsc10
           \fi
           \@scfont
}%
%
% The major parts of an entry are separated with \bblnewblock.  The
% numbers below are taken from LaTeX's `article' style.
%
\def\bblnewblock{\hskip .11em plus .33em minus .07em }%
%
% Here's where you stick any other bibliography-formatting goodies, or
% redefine the values above.
%
\let\bblhook = \empty
%
%
% Here are the four default definitions for formatting the in-text
% citations.  These are what you redefine (after your \input btxmac but
% before your \bibliography) to get parens instead of brackets, or
% superscripts, or footnotes, or whatever.
%
\def\printcitestart{[}%         left bracket
\def\printcitefinish{]}%        right bracket
\def\printbetweencitations{, }% comma, space
\def\printcitenote#1{, #1}%     comma, space, note (if it exists)
%
% That scheme is pretty flexible.  For example you could use
%     \def\printcitestart{\unskip $^\bgroup}
%     \def\printcitefinish{\egroup$}
%     \def\printbetweencitations{,}
%     \def\printcitenote#1{\hbox{\sevenrm\space (#1)}}
%     \font\eighttt = cmtt8
%     \scriptfont\ttfam = \eighttt
% to get superscripted in-text citations.  (The scriptfont stuff
% exists only to print an undefined citation; it's in cmtt8 because
% there is no cmtt7.)  To get something radically different, however,
% you'll have to define your own \cite command.
%
% When we read `\citation' from the .aux file, it means nothing.
%
\let\citation = \@gobble
%
%
% Now comes the stuff for dealing with LaTeX's \newcommand.  As
% mentioned earlier, this \newcommand will redefine a preexisting
% command; that's different from how LaTeX's \newcommand behaves.
%
\@innernewcount\@numparams
%
% \newcommand{\foo}[n]{text} defines the control sequence \foo to have
% n parameters, and replacement text `text'.
%
\def\newcommand#1{%
   \def\@commandname{#1}%
   \@getoptionalarg\@continuenewcommand
}%
%
% Figure out if this definition has parameters.
%
\def\@continuenewcommand{%
   % If no optional argument, we have zero parameters.  Otherwise, we
   % have that many.
   \@numparams = \ifx\@optionalarg\empty 0\else\@optionalarg \fi \relax
   \@newcommand
}%
%
% \@numparams is how many arguments this command has.  The name of the
% command is \@commandname.  The replacement text for the new macro is #1.
%
\def\@newcommand#1{%
   \def\@startdef{\expandafter\edef\@commandname}%
   \ifnum\@numparams=0
      \let\@paramdef = \empty
   \else
      \ifnum\@numparams>9
         \errmessage{\the\@numparams\space is too many parameters}%
      \else
         \ifnum\@numparams<0
            \errmessage{\the\@numparams\space is too few parameters}%
         \else
            \edef\@paramdef{%
               % This is disgusting, but \loop doesn't work inside \edef,
               % because \body isn't defined.
               \ifcase\@numparams
                  \empty  No arguments.
               \or ####1%
               \or ####1####2%
               \or ####1####2####3%
               \or ####1####2####3####4%
               \or ####1####2####3####4####5%
               \or ####1####2####3####4####5####6%
               \or ####1####2####3####4####5####6####7%
               \or ####1####2####3####4####5####6####7####8%
               \or ####1####2####3####4####5####6####7####8####9%
               \fi
            }%
         \fi
      \fi
   \fi
   \expandafter\@startdef\@paramdef{#1}%
}%
%
%% [[[end of BibTeX-specific stuff]]]
%
%
% Names of references (arguments given in the \cite and \nocite
% commands) and file names (arguments given in the \bibliography and
% \bibliographystyle commands) are recorded in \jobname.aux, called the
% \@auxfile in these macros.  Here's how they get read in.
%
\def\@readauxfile{%
   \if@auxfiledone \else % remember: \@auxfiledonetrue if \noauxfile is defined
      \global\@auxfiledonetrue
      \@testfileexistence{aux}%
      \if@fileexists
         \begingroup
            % Because we might be in horizontal mode when \@readauxfile
            % is called (if it's in response to a \cite or \nocite), we
            % want to ignore all the would-be spaces at the ends of
            % lines in the aux file.  Fortunately, it's highly unlikely
            % an end-of-line might actually be desired.
            % And because we don't change the category code of anything
            % but @, primitives like \gdef can't be used to define labels
            % in the aux file.  The solution adopted by btxmac.tex is to
            % write `\@citedef{LABEL}{DEFINITION}' to the aux file, and
            % use \csname on LABEL.
            \endlinechar = -1
            \catcode`@ = 11
            \input \jobname.aux
         \endgroup
      \else
         \message{\@undefinedmessage}%
         \global\@citewarningfalse
      \fi
      \immediate\openout\@auxfile = \jobname.aux
   \fi
}%
%
% The \@readauxfile macro does all that work the first time it's called.
% Since it's called once for every \cite, \nocite, \bibliography, and
% \bibliographystyle command that the user issues, we need to remember
% whether the work's been done.  It's considered done if we're not to do
% it---that is, if \noauxfile is defined.
%
\newif\if@auxfiledone
\ifx\noauxfile\@undefined \else \@auxfiledonetrue\fi
%
% It's conceivable you'd want to change how other characters are read;
% to do that, change their category code before doing \input btxmac.
%
%
% After reading the .aux file, \@readauxfile opens it for writing.
% The \@writeaux macro does the actual writing (as long as
% \noauxfile is undefined).
%
\@innernewwrite\@auxfile
\def\@writeaux#1{\ifx\noauxfile\@undefined \write\@auxfile{#1}\fi}%
%
%
% A macro package that uses btxmac.tex might define
% \@undefinedmessage (before doing an \input btxmac).
%
\ifx\@undefinedmessage\@undefined
   \def\@undefinedmessage{No .aux file; I won't give you warnings about
                          undefined citations.}%
\fi
%
% Even if citations are undefined, we want to complain only if
% \@citewarningtrue.  The default is to set \@citewarningtrue unless
% \noauxfile is defined.  Again, a macro package that uses
% btxmac.tex might want to redefine this.
%
\@innernewif\if@citewarning
\ifx\noauxfile\@undefined \@citewarningtrue\fi
%
%
% Finally, before leaving we restore @'s old category code.
%
\catcode`@ = \@oldatcatcode

\def\pfeilso{\leavevmode
            \vrule width 1pt height9pt depth 0pt\relax
           \vrule width 1pt height8.7pt depth 0pt\relax
           \vrule width 1pt height8.3pt depth 0pt\relax
           \vrule width 1pt height8.0pt depth 0pt\relax
           \vrule width 1pt height7.7pt depth 0pt\relax
            \vrule width 1pt height7.3pt depth 0pt\relax
            \vrule width 1pt height7.0pt depth 0pt\relax
            \vrule width 1pt height6.7pt depth 0pt\relax
            \vrule width 1pt height6.3pt depth 0pt\relax
            \vrule width 1pt height6.0pt depth 0pt\relax
            \vrule width 1pt height5.7pt depth 0pt\relax
            \vrule width 1pt height5.3pt depth 0pt\relax
            \vrule width 1pt height5.0pt depth 0pt\relax
            \vrule width 1pt height4.7pt depth 0pt\relax
            \vrule width 1pt height4.3pt depth 0pt\relax
            \vrule width 1pt height4.0pt depth 0pt\relax
            \vrule width 1pt height3.7pt depth 0pt\relax
            \vrule width 1pt height3.3pt depth 0pt\relax
            \vrule width 1pt height3.0pt depth 0pt\relax
            \vrule width 1pt height2.7pt depth 0pt\relax
            \vrule width 1pt height2.3pt depth 0pt\relax
            \vrule width 1pt height2.0pt depth 0pt\relax
            \vrule width 1pt height1.7pt depth 0pt\relax
            \vrule width 1pt height1.3pt depth 0pt\relax
            \vrule width 1pt height1.0pt depth 0pt\relax
            \vrule width 1pt height0.7pt depth 0pt\relax
            \vrule width 1pt height0.3pt depth 0pt\relax}

\def\pfeilsw{ \leavevmode 
            \vrule width 1pt height0.3pt depth 0pt\relax
            \vrule width 1pt height0.7pt depth 0pt\relax
            \vrule width 1pt height1.0pt depth 0pt\relax
            \vrule width 1pt height1.3pt depth 0pt\relax
            \vrule width 1pt height1.7pt depth 0pt\relax
            \vrule width 1pt height2.0pt depth 0pt\relax
            \vrule width 1pt height2.3pt depth 0pt\relax
            \vrule width 1pt height2.7pt depth 0pt\relax
            \vrule width 1pt height3.0pt depth 0pt\relax
            \vrule width 1pt height3.3pt depth 0pt\relax
            \vrule width 1pt height3.7pt depth 0pt\relax
            \vrule width 1pt height4.0pt depth 0pt\relax
            \vrule width 1pt height4.3pt depth 0pt\relax
            \vrule width 1pt height4.7pt depth 0pt\relax
            \vrule width 1pt height5.0pt depth 0pt\relax
            \vrule width 1pt height5.3pt depth 0pt\relax
            \vrule width 1pt height5.7pt depth 0pt\relax
            \vrule width 1pt height6.0pt depth 0pt\relax
            \vrule width 1pt height6.3pt depth 0pt\relax
            \vrule width 1pt height6.7pt depth 0pt\relax
            \vrule width 1pt height7.0pt depth 0pt\relax
            \vrule width 1pt height7.3pt depth 0pt\relax
            \vrule width 1pt height7.7pt depth 0pt\relax
            \vrule width 1pt height8.0pt depth 0pt\relax
            \vrule width 1pt height8.3pt depth 0pt\relax
            \vrule width 1pt height8.7pt depth 0pt\relax
            \vrule width 1pt height9pt depth 0pt\relax
      }

%  *** end including bib4plain.tex *** 
  % This will define \cite and make sure it works as in latex

\def\widestnumber#1#2{}
  % Our amstex-ppt style does not know about \widestnumber

\def\citewarning#1{\ifx\shlhetal\relax 
   % normal mode, do not write anything
    \else
   % private mode
    \par{#1}\par
    \fi
}

\def\rm{\fam0 \tenrm}

\def\fakesubhead#1\endsubhead{\bigskip\noindent{\bf#1}\par}

% % \input rsfs
%  *** start including rsfs.tex *** 

% # Keywords: Script or Calligraphic (Caligraphic) letters with the RSFS Font

% The story so far:    July 1998 -- Saharon would like to have a
% ``nicer'' calligraphic font. In particualr, the leters S and P in
% the usual calligraphic font do not look ``special'' enough. 
% 
% I found out that ``rsfs'' (``Ralph Smith Formal Script'') may be
% what he wants.   I installed the mf file, the .tfm file, as well as
% a few pk files in ~/TeX/rsfs.    Let's hope that this is enough.
% Using amstex, all you have to do is to \input rsfs.tex 
% Files prepared with citealice willdothis automatically. 
%
%  Note:  for some reason xdvi calls MakeTeXpk, then Maketexpk
%  complains about wrong resolution, but still writes commands to
%  missfont.log...  
%

% we redefine a macro inside amstex's \Cal command , so that it calls
% our nice font ``rsfs'' rather than the usual calligraphic font. 
% Note thisworks for amstex only.   
% In plain tex, would have to add definitions of \Cal
% in latex... we should insteaduse mathrsfs.sty
% 

\font\textrsfs=rsfs10
\font\scriptrsfs=rsfs7
\font\scriptscriptrsfs=rsfs5

\newfam\rsfsfam
\textfont\rsfsfam=\textrsfs
\scriptfont\rsfsfam=\scriptrsfs
\scriptscriptfont\rsfsfam=\scriptscriptrsfs

\edef\oldcatcodeofat{\the\catcode`\@}
\catcode`\@11

\def\Cal@@#1{\noaccents@ \fam \rsfsfam #1}

\catcode`\@\oldcatcodeofat

%  *** end including rsfs.tex *** 

\expandafter\ifx \csname margininit\endcsname \relax\else\margininit\fi

\long\def\red#1\endred{}
\long\def\green#1\endgreen{}
\long\def\blue#1\endblue{}
\long\def\private#1\endprivate{}

\def\endred{ \unmatched endred! }
\def\endgreen{ \unmatched endgreen! }
\def\endblue{ \unmatched endblue! }
\def\endprivate{ \unmatched endprivate! }

\ifx\latexcolors\undefinedcs\def\latexcolors{}\fi

\def\emptycs{}
\def\evaluatelatexcolors{%
        \ifx\latexcolors\emptycs\else
        \expandafter\xxevaluate\latexcolors\xxfertig\evaluatelatexcolors\fi}
\def\xxevaluate#1,#2\xxfertig{\setupthiscolor{#1}%
        \def\latexcolors{#2}}

 % \def\definedprivate{rgb 0.5 0 0.5} 
       % black 

\font\smallfont=cmsl7
\def\rutgerscolor{\ifmmode\else\endgraf\fi\smallfont% \vrule width 3cm height 1cm depth 0cm
\advance\leftskip0.5cm\relax}
\def\setupthiscolor#1{\edef\tmptmpcs{\noexpand\bgroup\noexpand\rutgerscolor
\noexpand\def\noexpand\currentcolor{#1}%
\noexpand}%
\expandafter\let\csname#1\endcsname\tmptmpcs
\def\tmptmpcs{\checkColorUnmatched{#1}\popthecolor}
\expandafter\let\csname end#1\endcsname\tmptmpcs}

\def\checkColorUnmatched#1{\def\expectcolor{#1}%
    \ifx\expectcolor\currentcolor   % OK! 
    \else \edef\failhere{\noexpand\tryingToClose '\currentcolor' with end\expectcolor}\failhere\fi}

\def\currentcolor{???}

\def\popthecolor{\ifmmode\else\endgraf\fi\egroup}

\expandafter\def\csname#1\endcsname{}

\evaluatelatexcolors

 \let\outerhead\head
 \def\head{\innerhead}
 \let\innerhead\outerhead

 \let\outersubhead\subhead
 \def\subhead{\innersubhead}
 \let\innersubhead\outersubhead

 \let\outersubsubhead\subsubhead
 \def\subsubhead{\innersubsubhead}
 \let\innersubsubhead\outersubsubhead

 \let\outerproclaim\proclaim
 \def\proclaim{\innerproclaim}
 \let\innerproclaim\outerproclaim

 % 
 % \newdimen\wzero 
 % \newdimen\hzero 
 % \newdimen\dzero 
 % 
 % \def\whiteblack#1{\ifx \shlhetal\relax      
 %    \hbox{#1}% normal mode
 %     \else\ifx\shlhetal\undefinedcontrolsequence
 %    \hbox{#1}% normal mode
 %         \else
 %    % private mode
 %     \setbox0=\hbox{#1}\leavevmode
 %     \wzero=\wd0\hzero=\ht0\dzero=\dp0
 %         \advance\wzero1mm
 %         \advance\dzero1mm
 %         \advance\hzero1mm
 %     \hbox{\vrule width \wzero height \hzero depth\dzero\hskip-\wzero
 %         \special{ps: currentgray 1  setgray}%
 %         \hbox to \wzero{\hss\copy0\hss}%
 %            \special{ps:setgray}}\fi\fi
 % }
 % 
 % 

\def\demo#1{\medskip\noindent{\it #1.\/}}
\def\enddemo{\smallskip}

\def\remark#1{\medskip\noindent{\it #1.\/}}
\def\endremark{\smallskip}

%  *** end including alice2jlem.tex *** 
\pageheight{8.5truein}
%\pageheight{48.5pc}
\topmatter
\title{Density is at most the spread of the square } \endtitle
\author {Saharon Shelah \thanks {\null\newline I would like to thank 
Alice Leonhardt for the beautiful typing. \null\newline
Publication E56} \endthanks} \endauthor 

\affil{The Hebrew University of Jerusalem \\
Einstein Institute of Mathematics \\
Edmond J. Safra Campus, Givat Ram \\
Jerusalem 91904, Israel
 \medskip
 Department of Mathematics \\
 Hill Center-Busch Campus \\
  Rutgers, The State University of New Jersey \\
 110 Frelinghuysen Road \\
 Piscataway, NJ 08854-8019 USA} \endaffil
\endtopmatter
\document

\newpage

\head {\S1} \endhead  \resetall \sectno=1
 \spuriousreset
\bigskip

\proclaim{\stag{1n.7} Claim}  Assume $\Bbb B$ is an infinite Boolean
Algebra and $\lambda = d(\Bbb B)$.  Then ${\frak s}(\Bbb B * \Bbb B)$, i.e.
${\frak s}(\text{\rm uf}(\Bbb B)
\times \text{\rm uf}(\Bbb B)) \ge \lambda$ (if $\lambda$ limit-obtained).
\endproclaim
\bigskip

\remark{Remark}  1) ul$(\Bbb B)$ is the space of ultrafilters of $\Bbb
B$ a compact space with clopen base.
\nl
2) ${\frak s}(X)$ is sup$\{|Y|:Y \subseteq X$ is discrete $Y$, same as
des$(X)\}$.
\nl
3) We meant to consider whether this works for compact Hausdorff
spaces.  
But subsequently and independently Szentmikl\"ossy prove this.
\endremark
\bigskip

\demo{Proof}  Without loss of generality $\lambda > \aleph_0$.  We
choose $(p^0_i,p^1_i,a_i)$ by induction on $i < \lambda$ such that
\mr
\item "{$\otimes$}"  $(a) \quad p^\ell_i$ is an ultrafilter of $\Bbb
  B$ for $\ell=0,1$
\sn
\item "{${{}}$}"  $(b) \quad a_i \in p^1_i,a_i \notin p^0_i$,
  i.e. $(-a_i) \in p^1_i$
\sn
\item "{${{}}$}"  $(c) \quad$ if $j<i$ then $a_j \in p^0_i
  \Leftrightarrow a_j \in p^1_i$
\sn 
\item "{${{}}$}"  $(d) \quad$ if $j<i$ then $a_i \notin p^0_j,a_i
  \notin p^1_j$.
\ermn
So assume we have arrived to $i$.  Let $\Bbb B_i$ be the subalgebra of
$\Bbb B$ generated by $\{a_j:j<i\}$.

For every non-zero $b \in \Bbb B_i$ choose an ultrafilter $q^i_b$ of
$\Bbb B$ and for simplicity $b = a_j \Rightarrow q^i_b = p^1_j$ and $b
= (-a_j) \Rightarrow q^i_0 = p^0_j$ for $j<i$.

As $d(\Bbb B) \ge \lambda$ clearly $\{q^i_b:b \in \Bbb B_i \backslash
\{0\}\}$ is not dense hence there is a non-zero $a_i \in \Bbb B$ such
that $b \in \Bbb B_i \backslash \{0\} \Rightarrow a_i \notin q^i_b$
(i.e. a non-empty clopen set to which none of the points $q^i_p$
belongs).

Now clearly $b \in \Bbb B_i \backslash \{0\} \Rightarrow a_i \ne b$
(as $b \in q^i_p$) hence $a_i \notin \Bbb B_i$.  This implies that
there is an ultrafilter $q^*_i$ of $\Bbb B_i$ such that
\mr 
\item "{$\circledast$}"  $b \in q \Rightarrow a_i \cap b > 0 \wedge
(-a_i) \cap b >0$.
\ermn
[Why?  As $\{b_0 \cup b_1:b_0,b_1 \in 
\Bbb B_2 \text{ and } b-1 \cap a_i = 0_{\Bbb B}\text{ and } b_a - a_i = 0\}$
is a proper ideal of $\Bbb B_i$ hence can be extended to an
ultrafilter of $\Bbb B_i$.]

So there are ultrafilters $p^0_i,p^1_i$ of $\Bbb B$ such that
\mr 
\item "{$\circledast$}"  $q^*_i \cup \{a_i\} \subseteq p^1_i$ and
  $q^*_i \cup \{-a_i\} \subseteq p^0_i$.
\ermn
This is enough for the induction step.

Having carried the induction
\mr
\item "{$(a)$}"  $p_i := (p^0_i,p^1_i) \in \text{ uf}(\Bbb B) \times
  \text{ uf}(\Bbb B)$
\sn 
\item "{$(b)$}"  $(-a_i) \times a_i$ is an open subset of uf$(\Bbb B)
  \times \text{ uf}(\Bbb B)$.
\ermn
Lastly,
\mr
\item "{$(c)$}"  if $i < j < \lambda$ \ub{then} $p_j \notin (-a_i)
\times a_i$ because $(-a_i) \notin p^0_j$ or $a_i \notin p^1_j$ as
$a_i \in \Bbb B_j$ and $p^0_j \cap B_j = p^1_j \cap \Bbb B_j$ by the
choice of $p^0_j,p^1_0$
\sn 
\item "{$(d)$}"  if $i< j < \lambda$ then $p_i \notin (-a_j) \times
  a_j$ because $p^1_i \notin a_j$ by the choice of $a_j$.
\ermn
So $\langle(p_i,(-a_i) \times a_i):i < \lambda)\rangle$ exemplify
dis$(\text{uf}(\Bbb B) \times \text{ uf}(\Bbb B)) \ge \lambda$ as required.
\enddemo
\newpage

% PRIVATE excursion starts here
% back from PRIVATE part number 1/1

\enddocument